\documentclass[11pt]{article} 
\usepackage{amsmath}
\usepackage{amssymb}     
\usepackage{latexsym}   

\input diagrams.tex

\def\endproof{$\ \ \Box$} 

\newcommand{\bq}{\begin{quote}}
\newcommand{\eq}{\end{quote}}

\newcommand{\dom}{\mathrm{dom}}
\newcommand{\fix}{\mathrm{fix}}

\newcommand{\reals}{{\mathbb R}}

\newtheorem{Th}{Theorem}[section]
\newtheorem{theorem}[Th]{Theorem}
\newtheorem{proposition}[Th]{Proposition}     
\newtheorem{lemma}[Th]{Lemma}

\newtheorem{definition}[Th]{Definition}

\newtheorem{remark}[Th]{Remark}

\begin{document}
\title{The maximum entropy state}
\author{Keye Martin \\ \\
{\small Department of Mathematics}\\
{\small Tulane University}\\
{\small New Orleans, Louisiana 70118}\\
{\small \texttt{martin@math.tulane.edu } }\\
{\small \texttt{http://www.math.tulane.edu/$\tilde{\ }$martin}}}
\date{}
\maketitle
\begin{abstract} We give 
an algorithm for calculating
the maximum entropy state as 
the least fixed point of a Scott continuous 
mapping on the domain of classical states 
in their Bayesian order.
\end{abstract}

\section{Introduction}

Suppose an experiment has one of $n$
different possible outcomes. These outcomes define a function $a:\{1,\ldots,n\}\rightarrow\reals$
we sometimes call an \em observable\em.
Now suppose in addition that we know
that if this experiment were performed repeatedly
that the average outcome $\langle a\rangle$ would be $E$. 
For any number of
reasons we can imagine wanting to determine a distribution
$x\in\Delta^n$ with $\langle a|x\rangle=E$ that
is a good candidate for the `actual probabilities.' 
As an example, if we have a six sided die, and we know
the average roll is 3.5, then we could all agree
that the `best distribution' which matches the information
we have is $(1/6,\ldots,1/6)$, i.e., all sides are
equally likely. But now suppose we know the die
is biased in some way and that the average $\langle a\rangle$
is $E\not =3.5$?

The
difficulty mathematically is that we only have two equations
($x\in\Delta^n$ and $\langle a|x\rangle=E$) but that we are trying to
solve for $n$ unknowns. The \em maximum entropy principle \em
offers an approach to this problem: it says we should
choose the unique $x\in\Delta^n$ with $\langle a|x\rangle=E$
whose entropy $\sigma x$ is maximum. Later 
we will give a proof that the problem has a unique
solution and how to find it (because a complete
detailed proof
of this well known result either does not exist
in the literature or is very good at hiding).

The maximum entropy principle originates from
statistical mechanics, where it was
known at least as far back as Gibbs that
the Boltzmann distribution is the unique
distribution which maximizes
the entropy while simultaneously
also yielding the observed average energy ($\propto$ temperature). 
But the development
which appears to have led Jaynes in \cite{jaynes:I} 
to regard the maximum entropy principle as
a legitimate form of inference was Shannon's
information theory~\cite{shannon} and its
characterization of entropy as the
unique function satisfying three axioms
that a measure of uncertainty might obey.
In~\cite{jaynes:I}, Jaynes offered a
very different explanation of statistical
mechanics: it has a physical part,
which serves to enumerate the states
of a system and their properties,
and it has an \em inferential \em
aspect, where conclusions are drawn
on the basis of incomplete information by 
(a) accepting entropy as the canonical
measure of uncertainty in a distribution, partly
because of Shannon's theorem, and  
(b) applying the maximum entropy principle.
A benefit of this viewpoint is that the
nonphysical assumptions required by traditional approaches to
the subject -- assumptions not derivable from the laws of motion --
are replaced by a simple inference principle 
and we can still derive
the usual computational rules in statistical
mechanics, such as the Boltzmann distribution~\cite{jaynes:I}.
Since Jaynes'
endorsement of the idea~\cite{jaynes:I}, 
the maximum entropy principle has been
successfully applied to problems in many disciplines, including
spectral analysis~\cite{burg}, image reconstruction (\cite{im1}\cite{im2}),
and somewhat recently to natural language processing (\cite{natural1}\cite{natural2}).

There are many instances in the literature
where it is necessary to calculate the maximum
entropy state, and the methods we have seen
look to be rather geared toward the application
at hand. It is known that calculating this state
amounts to being able to solve the equation
\[E\cdot\sum_{i=1}^n e^{\lambda a_i}=\sum_{i=1}^n a_i e^{\lambda a_i}\]
for the \em Lagrangian multiplier \em $\lambda$.
 But we have not seen any
methods which have been proven to always work.
(For instance, one could ask, does Newton's
method work, and if so, for which initial guesses?
We don't know the answer to this question by the
way, but we suspect that it always does.) In this paper 
we give a method for
calculating the maximum entropy state that
always applies, we prove that it works starting
from any initial guess, and that it is really
a domain theoretic idea in disguise:
the maximum entropy state is the least fixed point 
of a Scott continuous map on the domain of classical 
states in their Bayesian order.

The heart of the method makes use of
basic techniques from numerical analysis.
What requires work is to \em realize \em that these
basic techniques are perfectly suited for
solving this nontrivial equation. This realization
is only possible if one is persistent enough 
to continue manipulating
some really dreadful sums until they take
just the right form. All the arithmetic
pays off in the end because there are some
really neat domain theoretic ideas hidden 
beneath the surface useful for calculating
the Lagrangian multiplier that defines
the maximum entropy state.

\section{The equilibrium state}

For an integer $n\geq 2$, the \em classical states \em are 
\[\Delta^n:=\left\{x\in[0,1]^n:\sum_{i=1}^nx_i=1\right\}.\] 
A classical state $x\in\Delta^n$ is \em pure \em when $x_i=1$ for some 
$i\in\{1,\ldots,n\}$; we denote such a state by $e_i$. 
Pure states $\{e_i\}_i$ are the actual states a system can 
be in, while general mixed states $x$ and $y$ are epistemic entities. 

Given a vector $a:\{1,\ldots,n\}\rightarrow\reals$, sometimes
called an \em observable\em, its average value is
\[\langle a|x\rangle = \sum_{i=1}^n a_i x_i\]
defined for $x\in\Delta^n$. Shannon entropy
$\sigma:\Delta^n\rightarrow[0,\infty)$ is
\[\sigma x=-\sum_{i=1}^n x_i\log x_i.\]
These ideas combine
to give the existence
and uniqueness of the equilibrium state associated
to (energy) observable $a$ in thermodynamics.

\begin{lemma} 
\label{eq} If $a:\{1,\ldots,n\}\rightarrow\reals$ is a vector, 
there is a unique classical state $y\in\Delta^n$ such that
\[\langle a|y\rangle-\sigma y=\inf\{\langle a|x\rangle-\sigma x:x\in\Delta^n\}.\]
The state $y$ is given pointwise by $y_i=e^{-a_i}/Za$ and satisfies
\[\langle a|y\rangle-\sigma y=-\log Za\]
where
\[Za:=\sum_{i=1}^n\frac{1}{e^{a_i}}.\]
\end{lemma}
{\bf Proof}. First, arithmetic gives $\langle a|y\rangle-\sigma y=-\log Za.$
Next, it is the minimum value of $f(x)=\langle a|x\rangle-\sigma x$
on $\Delta^n$:
\begin{eqnarray*}
f(x) & = & f(x) + \log Za-\log Za\\
     & = & f(x) + \sum_{i=1}^n (x_i\log Za)-\log Za\\
     & = & -\sum_{x_i>0}\log\left(\frac{y_i}{x_i}\right)x_i -\log Za\\
     & \geq &\sum_{x_i>0}\left(1-\frac{y_i}{x_i}\right)x_i-\log Za\ \
     \  
     (\mbox{using}\ \log x\leq x-1\ \mbox{for}\ x>0)\\
     & = & \left(1-\sum_{x_i>0}y_i\right)-\log Za\\
     &\geq& -\log Za
\end{eqnarray*}
Finally, $y$ is the \em unique state \em where $f$ takes it
minimum: If $f(x)=-\log Za$, then the string of inequalities
above implies
\[-\sum_{x_i>0}\log\left(\frac{y_i}{x_i}\right)x_i=\sum_{x_i>0}\left(1-\frac{y_i}{x_i}\right)x_i\]
which can be rewritten as
\[\sum_{x_i>0}(t_i-1-\log t_i)x_i=0\]
where $t_i=y_i/x_i$. Because $\log x\leq x-1$ for $x>0$, 
this is a sum of nonnegative terms which results in zero. Then
each term must be zero, so $t_i=1$ which means $x_i=y_i$ whenever $x_i>0$.
However, since $\sum y_i=1$ and each $y_i>0$,
we must have $x_i>0$ for \em all \em $i\in\{1,\ldots,n\}.$ Then $x=y$.
\endproof\newline

\section{The maximum entropy principle}

For the rest of the paper, we now
fix an observable $a:\{1,\ldots,n\}\rightarrow\reals$ with $a_i<a_{i+1}$.
If we know the probability $x_i$ that outcome $a_i$ will occur,
then we have a distribution $x\in\Delta^n$ which
can be used to calculate the average value $\langle a\rangle$ of $a$:
\[\langle a\rangle=\langle a|x\rangle = \sum_{i=1}^n a_i x_i.\]
But suppose all we know is that $\langle a\rangle=E$. What
distribution should we attribute
this average to?  The maximum entropy principle says 
we should choose the $x\in\Delta^n$ with $\langle a|x\rangle=E$
whose entropy $\sigma x$
is as large as possible. Then we seek
to maximize $\sigma x$
subject to the constraints
\[\sum_{i=1}^n x_i=1\ \ \ \mbox{and}\ \ \ \langle a|x\rangle = \sum_{i=1}^n a_i x_i = E.\]
Let us think about this carefully. First,
assuming for the moment that $S=\langle a|\cdot\rangle^{-1}(\{E\})\not=\emptyset$,
the problem has a solution because $\sigma$
is continuous on a nonempty closed subset 
$S$ of the compact set $\Delta^n$. Call
this solution $x$. If there were another
solution $y\neq x$, then it too
would satisfy $\langle a|y\rangle=E$ and $\sigma y=\sigma x$.
But then $z=(x+y)/2\in\Delta^n$
would have $\langle a|z\rangle = E$ and 
by the strict concavity of entropy,
\[\sigma(z)> (1/2)\sigma x+(1/2)\sigma y=\sigma x,\]
which contradicts the fact that $x$ maximizes
the entropy. So the problem
has a \em unique \em solution
and we call it \em the maximum entropy state. \em

In every reference we have seen in the literature,
the next step has been to apply Lagrangian multipliers
to determine a candidate for the maximum of $\sigma$ on $S$. 
But Lagrangian multipliers only applies to regions 
where all partial derivatives of $\sigma$ (and the
two constraints) exist,
and then it is only capable of  detecting extrema that occur
at interior points of such a region. 
The partial derivatives of $\sigma$ do not exist on the boundary of $\Delta^n$.
So if we want a guarantee
that Lagrangian multipliers will yield the maximum entropy
state, we need to know that the maximum entropy state does not occur
on the boundary of $\Delta^n$. From the
point of view of a pure optimization problem,
it is not clear why the maximum should always be taken
on the interior of $\Delta^n$. Though
it proves to be largely technical, we still
should make the following point: it is mathematically
incorrect to apply Lagrangian multipliers in
this situation and then assert that it always
 yields the maximum entropy state. 

As an example, if $E=a_1$, then $e_1$ is
the maximum entropy state, which lies along
the boundary. For $E=a_n$, the maximum entropy state is $e_n$. Ignoring
this for a moment, a suspect 
application of Lagrangian multipliers suggests 
that the maximum entropy state $y$
is given by
\[y_i=\frac{e^{\lambda a_i}}{\sum_{i=1}^n e^{\lambda a_i}}\]
where $\lambda\in\reals$ satisfies
\[E\cdot\sum_{i=1}^n e^{\lambda a_i}=\sum_{i=1}^n a_i e^{\lambda a_i}.\]
At this stage, we have 
no way of knowing which one of the following is true:
\begin{enumerate}
\item[(i)] The maximum entropy state occurs on 
the boundary, in which case it is not the state $y$,
 so Lagrangian multipliers
is of no use in finding it,
\item[(ii)] The maximum entropy state occurs
at an interior point, in which case the equation
for $\lambda$ has at least one solution, and one
of these solutions will yield the maximum entropy state. 
But which one?
\end{enumerate}

Though $\lambda$ may not even exist, it turns out that there are only two cases where this is true:
$E=a_1$ and $E=a_n$. Otherwise, $\lambda$ exists uniquely,
and defines the unique state $y$, leaving three possibilities: $y$ is the
maximum entropy state, a minimum, or neither.  Thankfully:

\begin{proposition} Let $a:\{1,\ldots,n\}\rightarrow\reals$ be a vector with $a_i<a_{i+1}$
and $a_1<E<a_n$. There is a unique $\lambda\in\reals$ with
\[E\cdot\sum_{i=1}^n e^{\lambda a_i}=\sum_{i=1}^n a_i e^{\lambda a_i}.\]
The state $y\in\Delta^n$ given by
\[y_i=\frac{e^{\lambda a_i}}{\sum_{i=1}^n e^{\lambda a_i}}\]
satisfies 
\[\langle a|y\rangle=E\ \ \mbox{and}\ \ \sigma y = \sup\{\sigma x:\langle a|x\rangle=E\ \&\ x\in\Delta^n\}\]
Thus, $y$ is
the only state with these two properties, i.e.,
it is the maximum entropy state.
\end{proposition}
{\bf Proof}. First suppose that a solution $\lambda$ to the equation
exists. We will prove that the associated $y$ has
to be the maximum entropy state. To do so, define a new
observable $b=-\lambda a$ by $b_i=-\lambda a_i$.
Now take the equilibrium state $y$ associated to $b$ given by
\[y_i=\frac{e^{\lambda a_i}}{Zb}\ \ \ \&\ \ \ Zb=\sum_{i=1}^n e^{\lambda a_i}.\]
By Lemma~\ref{eq} we also have
\[\langle b|y\rangle-\sigma y =\inf\{\langle b|x\rangle-\sigma x:x\in\Delta^n\}.\]
Now we can prove that $y$ is the maximum entropy state.
First,
\[\langle a|y\rangle=\sum_{i=1}^n a_i\left(\frac{e^{\lambda a_i}}{Zb}\right)=\frac{\sum_{i=1}^n a_i e^{\lambda a_i}}{Zb}=E\]
using our assumption about $\lambda$. Next, let $x$ be any other
state with $\langle a|x\rangle=E$. Then $\langle b|x\rangle=-\lambda E = \langle b|y\rangle$.
Thus,
\[-\lambda E-\sigma x = \langle b|x\rangle - \sigma x \geq \langle b|y\rangle - \sigma y = -\lambda E - \sigma y\]
because $y$ is the equilibrium state for $b$. This
proves $\sigma y\geq\sigma x$.

Now suppose the equation has two solutions $\lambda$ and $\beta$.
Then each defines the maximum entropy state (which we know is unique) so
for all $i$,
\[y_i=\frac{e^{\lambda a_i}}{\sum_{i=1}^n e^{\lambda a_i}}=\frac{e^{\beta a_i}}{\sum_{i=1}^n e^{\beta a_i}}\]
In particular,
\[\frac{y_1}{y_2}=e^{\lambda(a_1-a_2)}=e^{\beta(a_1-a_2)}\]
which means $\lambda(a_1-a_2)=\beta(a_1-a_2)$, 
and since $a_1<a_2$, we get $\lambda=\beta$. Thus, $\lambda$
is unique assuming it exists.

Last, we prove that a solution to the equation exists. Define
\[f(x)=\frac{\sum_{i=1}^n a_i e^{x a_i}}{\sum_{i=1}^n e^{xa_i}} - E.\]
Using $a_1<a_n$, we have
\[\lim_{x\rightarrow -\infty} f(x)=a_1-E\ \ \ \&\ \ \ \lim_{x\rightarrow \infty} f(x)=a_n-E.\]
Because $a_1<E<a_n$, taking $c$ close to $-\infty$ gives $f(c)<0$
while $d$ close to $\infty$ gives $f(d)>0$. 
The continuity of $f$ yields a $\lambda$ with $f(\lambda)=0$.
\endproof\newline

For the case $E=a_1$, 
the only state $x$ with $\langle a|x\rangle=a_1$ is $x=e_1$,
so the maximum entropy state is $e_1$; for $E=a_n$, it is $e_n$. Then
\[(\exists x)\,\langle a|x\rangle=E\ \Leftrightarrow\ a_1\leq E\leq a_n\]
So, \em the maximum entropy state exists
if and only if $E\in[a_1,a_n]$. \em This is very
intuitive because it says that there is a solution
iff the expectation lies between the smallest and 
largest observable values. Going through
a proof of this well known result in detail
provides us with some of the technical
ideas  that will be useful in designing an algorithm
for actually calculating the maximum entropy state.

\section{An algorithm for calculating $\lambda$}
Define
\[f(x)=\frac{\sum_{i=1}^n a_i e^{x a_i}}{\sum_{i=1}^n e^{xa_i}} - E\]
and
\[I_f(x) = x - \frac{f(x)}{(a_n-a_1)^2}\]
for any $x\in\reals$.

\begin{lemma} Let $a_1<E<a_n$.
\begin{enumerate}
\em\item[(i)]\em For all $x\in\reals$, we have $0<I_f^\prime(x)<1$.
\em\item[(ii)]\em The map $I_f$ has a unique fixed point $\lambda$ and $I_f^n(x)\rightarrow \lambda$
for each $x\in\reals$.
\end{enumerate}
\end{lemma}
{\bf Proof}. (i) To prove $0<f^\prime(x)<(a_n-a_1)^2$, we first calculate $f^\prime(x)$. This
takes a while if we want to simply it as much as possible:
\begin{eqnarray*}
f^\prime(x) & = & \frac{1}{(\sum_i e^{xa_i})^2}\cdot\left(\sum_{i=1}^n e^{x a_i}\sum_{i=1}^n a_i^2 e^{xa_i}-\sum_{i=1}^n a_i e^{xa_i}
\sum_{i=1}^n a_i e^{x a_i}\right)\\
            & = & \frac{1}{(\sum_i e^{xa_i})^2}\cdot\left(\sum_{1\leq i,j\leq n} e^{xa_i}e^{xa_j}a_j^2-\sum_{1\leq i,j\leq n}a_ia_je^{xa_i}e^{xa_j}\right)\\
            & = & \frac{1}{(\sum_i e^{xa_i})^2}\cdot\left(\sum_{1\leq i,j\leq n}e^{xa_i}e^{xa_j}(a_j^2-a_ia_j)\right)\\
            & = & \frac{1}{(\sum_i e^{xa_i})^2}\cdot\left(\sum_{1\leq i\neq j\leq n}e^{xa_i}e^{xa_j}(a_j^2-a_ia_j)\right)\\
            & = & \frac{1}{(\sum_i e^{xa_i})^2}\cdot\left(\sum_{1\leq i< j\leq n}e^{xa_i}e^{xa_j}(a_j^2-a_ia_j)+\sum_{1\leq j<i\leq n}e^{xa_i}e^{xa_j}(a_j^2-a_ia_j) \right)\\
            & = & \frac{1}{(\sum_i e^{xa_i})^2}\cdot\left(\sum_{1\leq i< j\leq n}e^{xa_i}e^{xa_j}(a_j^2-a_ia_j)+\sum_{1\leq i<j\leq n}e^{xa_j}e^{xa_i}(a_i^2-a_ja_i) \right)\\
            & = & \frac{1}{(\sum_i e^{xa_i})^2}\cdot\left(\sum_{1\leq i< j\leq n}e^{xa_i}e^{xa_j}(a_j-a_i)^2\right)\\
\end{eqnarray*}
This proves $f^\prime(x)>0$. To prove $f^\prime(x)<(a_n-a_1)^2$,
\begin{eqnarray*}
f^\prime(x) & = & \frac{1}{(\sum_i e^{xa_i})^2}\cdot\left(\sum_{1\leq i< j\leq n}e^{xa_i}e^{xa_j}(a_j-a_i)^2\right)\\
            & < &  \frac{(a_n-a_1)^2}{(\sum_i e^{xa_i})^2}\cdot\left(\sum_{1\leq i< j\leq n}e^{xa_i}e^{xa_j}\right)\\
            & < &  (a_n-a_1)^2\\
\end{eqnarray*}
where the first inequality follows from the increasingness of $a$,
and the second inequality uses
\[\left(\sum_{i=1}^n e^{xa_i}\right)^2=\sum_{1\leq i,j\leq n} e^{xa_i}e^{xa_j} > \sum_{1\leq i< j\leq n}e^{xa_i}e^{xa_j}.\]
Then $0<f^\prime(x)<(a_n-a_1)^2$ which gives $0<I_f^\prime(x)<1$.\newline

(ii) For the $\lambda$ with $f(\lambda)=0$ we have $I_f(\lambda)=\lambda$. 
Given $x\in\reals$, set
\[I_{x}=[\min\{x,\lambda\},\max\{x,\lambda\}]\]
\[c_{x}=\sup_{t\in I_{x}} I_f^\prime(t).\]
Because $I_f^\prime$ is a continuous function
on a compact set $I_{x}$, it assumes its
absolute maximum at some point $t^*\in I_{x}$, i.e.,
$c_{x}=I_f^\prime(t^*)$. This proves $0<c_{x}<1$.
Then by the mean value theorem and
the fact that $I_f^\prime>0$,
\[d(I_f(a),I_f(b))\leq c_{x} \cdot d(a,b)\]
for all $a,b\in I_{x}$,
where $d$ is the usual metric on $\reals$. 
But $I_f(I_{x})\subseteq I_{x}$ because
$I_f$ maps sets of the form $[x,\lambda]$ and $[\lambda,x]$ to themselves,
using the strict monotonicity of $f$ that follows from $f^\prime>0$,
and the two equivalences
\[x\leq I_f(x)\equiv x\leq\lambda\ \ \ \&\ \ \ I_f(x)\leq x\equiv \lambda\leq x.\]
Thus, $I_f$ is a contraction on $I_x$, so it
has a unique fixed point on $I_x$, which must
be $\lambda$, and $I_f^n(x)\rightarrow\lambda$.
\endproof\newline

\begin{remark}\em $I_f$ is not a contraction because 
its derivative gets arbitrarily
close to one:
\[\lim_{x\rightarrow -\infty}I_f^\prime(x)=\lim_{x\rightarrow \infty}I_f^\prime(x)=1\]
Any contraction constant $c<1$
would have to bound its derivative from above. A forthcoming
work will study why `functions like these' have canonical fixed points.
\end{remark}

The equivalences
\[x\leq I_f(x)\equiv x\leq\lambda\ \ \ \&\ \ \ I_f(x)\leq x\equiv \lambda\leq x.\]
are important because they allow
us to determine properties of $\lambda$ without
actually knowing $\lambda$. For instance,
$\lambda >0$ iff $I_f(0)>0$. The advantage
is that $I_f(0)>0$ can be determined computationally,
while testing $\lambda>0$ would require us to
know the value of $\lambda$.

\section{From bottom to the maximum entropy state}

The calculation of $\lambda$ given in
the last section also has a formulation in terms of
classical states: the maximum entropy state is the
least fixed point of a Scott continuous map on the
$\Delta^n$ in its Bayesian order.

\begin{definition}\em A \em poset \em $P$
is a partially ordered set. A nonempty
subset $S\subseteq P$ is \em directed \em if
$(\forall x,y\in S)(\exists z\in S)\,x,y\sqsubseteq z$. The
\em supremum \em $\bigsqcup S$ of $S\subseteq P$ is
the least of its upper bounds when it exists.
A \em dcpo \em is a poset in which every directed set has a supremum.
\end{definition}
A function $f:D\rightarrow E$ between
dcpo's is \em Scott continuous \em if it is \em monotone\em,
\[(\forall x,y\in D)\,x\sqsubseteq y\Rightarrow f(x)\sqsubseteq f(y),\]
and \em preserves directed suprema\em,
\[f(\bigsqcup S)=\bigsqcup f(S),\]
for all directed $S\subseteq D$. Like complete metric spaces,
dcpo's also have a result which guarantees the existence
of \em canonical \em fixed points.

\begin{theorem} 
\label{leastfixedpoint}
Let $D$ be a dcpo with a least element $\bot$.
If $f:D\rightarrow D$ is a Scott continuous map, then
\[\mathrm{fix}(f):=\bigsqcup_{n\geq 0}f^n(\bot)\]
is the least fixed point of $f$ on $D$.
\end{theorem}

The set of classical states $\Delta^n$
has a natural domain theoretic structure,
too many of them in fact. The one
of interest in this paper is the Bayesian
order~\cite{meandbob}, which we now
briefly consider.
 
 Imagine that one of $n$ different outcomes is possible.
 If our knowledge of the outcome is $x\in\Delta^n$,
and then by some means we determine that outcome $i$ is 
not possible, our knowledge 
improves to 
\[p_i(x)=\frac{1}{1-x_i}(x_1,\ldots,\hat{x_i},\ldots,x_{n+1})\in\Delta^n,\] 
where $p_i(x)$ is obtained 
by first removing $x_i$ from $x$ and 
then renormalizing. The partial mappings which result, ${p_i:\Delta^{n+1}\rightharpoonup\Delta^n}$ 
with ${\dom(p_i)=\Delta^{n+1}\setminus\{e_i\}}$, 
are called the {\it Bayesian projections\,} and lead 
one to the following relation on classical states. 

\begin{definition}\em For $x,y\in\Delta^{n+1}$, 
\[x\sqsubseteq y\equiv(\forall 
i)(x,y\in\mbox{dom}(p_i)\Rightarrow p_i(x)\sqsubseteq p_i(y)).\]  
For $x,y\in\Delta^2$,  
\[x\sqsubseteq y\equiv (y_1\leq x_1\leq 1/2)\mbox{ or }(1/2\leq x_1\leq
y_1)\,.\] 
The relation $\sqsubseteq$ on $\Delta^n$ is called the \em Bayesian 
order. \em 
\end{definition} 
The Bayesian order was invented in~\cite{meandbob}
where the following is proven:

\begin{theorem} $(\Delta^n,\sqsubseteq)$ is a dcpo
with least element $\bot:=(1/n,\ldots,1/n)$ and $\mathrm{max}(\Delta^n)=\{e_i:1\leq i\leq n\}.$
\end{theorem}

The Bayesian order has a more direct description: The \em symmetric
formulation \em\cite{meandbob}. 
Let $S(n)$ denote the group of permutations on $\{1,\ldots,n\}$ 
and 
\[{\Lambda^n:=\{x\in\Delta^n:(\forall i<n)\,x_i\geq x_{i+1}\}}\] 
denote the collection of \em monotone \em decreasing
classical states. 

\begin{theorem}
\label{classicalsymmetries} 
For $x,y\in\Delta^n$, we have $x\sqsubseteq y$ iff 
there is a permutation ${\sigma\in S(n)}$ such that 
$x\cdot\sigma,y\cdot\sigma\in\Lambda^n$  
and 
\[(x\cdot\sigma)_i(y\cdot\sigma)_{i+1}\leq
(x\cdot\sigma)_{i+1}(y\cdot\sigma)_i\] 
for all $i$ with $1\leq i<n$. 
\end{theorem} 
Thus, $(\Delta^n,\sqsubseteq)$
can be thought of as $n!$ many copies of the domain $(\Lambda^n,\sqsubseteq)$
identified along their common boundaries, where $(\Lambda^n,\sqsubseteq)$ is
\[x\sqsubseteq y\equiv(\forall i<n)\,x_iy_{i+1}\leq x_{i+1}y_i.\]
It should be remarked though
that the problems of ordering $\Lambda^n$ and $\Delta^n$ are
very different, with the latter being far more challenging,
especially if one also wants to consider quantum mixed states.
Now to the fixed point theorem.

\begin{definition}\em Define $\lambda:\Delta^n\rightarrow\reals\cup\{\pm\infty\}$ by
\[\lambda(x)=\left\{ \begin{array}{ll}
                     \frac{\log\left(\frac{\mathrm{sort}(x)_1}{\mathrm{sort}(x)_2}\right)}{a_n-a_{n-1}}\ \mbox{\ if\ } I_f(0)>0;\\ \\
                     \frac{\log\left(\frac{\mathrm{sort}(x)_1}{\mathrm{sort}(x)_2}\right)}{a_1-a_{2}} \ \mbox{\ otherwise}.
                     \end{array}\right. \]
with the understanding for pure states that $\lambda x=\infty$ in the first
case and $\lambda x=-\infty$ in the other. The map $\mathrm{sort}$ puts
states into \em decreasing \em order.
\end{definition}

\begin{lemma} 
\label{lambdalemma}
For a function $a:\{1,\ldots,n\}\rightarrow\reals$ with $a_i<a_{i+1}$,
\begin{enumerate}
\em\item[(i)]\em If $\lambda>0$, then
\[(\forall x,y\in\Delta^n)\,x\sqsubseteq y\Rightarrow \lambda x\leq\lambda y\]
in the Bayesian order on $\Delta^n$.
\em\item[(ii)]\em If $\lambda\leq 0$, then
\[(\forall x,y\in\Delta^n)\,x\sqsubseteq y\Rightarrow \lambda x\geq\lambda y\]
in the Bayesian order on $\Delta^n$.
\end{enumerate}
\end{lemma}
           That is, the sign of the fixed point $\lambda=I_f(\lambda)$
           determines whether $\lambda:\Delta^n\rightarrow\reals$
           is monotone increasing ($\lambda >0$) or monotone decreasing ($\lambda\leq 0$). 
               
\begin{theorem} Let $a_1<E<a_n$. The map
$\phi:\Delta^n\rightarrow\Delta^n$ given by
\[\phi(x)=(e^{I_f(\lambda x)a_1},\ldots,e^{I_f(\lambda x)a_n})\cdot\frac{1}{Z(x)}\]
\[Z(x)=\sum_{i=1}^n e^{I_f(\lambda x)a_i}\]
is Scott continuous in the Bayesian order. Its least fixed point
is the maximum entropy state.
\end{theorem}      
{\bf Proof}. Let $x\sqsubseteq y$ in the Bayesian order. 
If $y\in\max(\Delta^n)$, then 
since $\lambda y=\pm\infty$,
the intent of the definition
is a limit. The state $\phi y$
is either $e_n$ or $e_1$. It follows that $\phi x\sqsubseteq\phi y$. Assume
now that $y\not\in\max(\Delta^n)$.

Now suppose $\lambda>0$. Then $0=\lambda\bot\leq\lambda x\leq\lambda y$
so $0<I_f(0)\leq I_f(\lambda x)\leq I_f(\lambda y)$. Then $\phi(x)$ and $\phi(y)$
are increasing states and we get
\begin{eqnarray*}
\phi(x)\sqsubseteq\phi(y) & \Leftrightarrow & (\forall 1\leq i< n)\, I_f(\lambda x)(a_{i+1}-a_i)\leq I_f(\lambda y)(a_{i+1}-a_i) \\
                          & \Leftrightarrow & I_f(\lambda x)\leq I_f(\lambda y)\\
                          & \Leftrightarrow  & \lambda x\leq \lambda y
\end{eqnarray*}
and this is true since $\lambda>0$ implies
$\lambda:\Delta^n\rightarrow\reals$ is monotone increasing. For $\lambda\leq 0$,
$\lambda y\leq\lambda x\leq\lambda\bot=0$, 
so $I_f(\lambda y)\leq I_f(\lambda x)\leq I_f(0)\leq 0$.
Then $\phi(x)$ and $\phi(y)$ are decreasing states and we get
\begin{eqnarray*}
\phi(x)\sqsubseteq\phi(y) & \Leftrightarrow & (\forall 1\leq i<n)\, I_f(\lambda x)(a_i-a_{i+1})\leq I_f(\lambda y)(a_{i}-a_{i+1}) \\
                          & \Leftrightarrow & I_f(\lambda x)\geq I_f(\lambda y)\\
                          & \Leftrightarrow  & \lambda x\geq \lambda y
\end{eqnarray*}
and this is true since $\lambda\leq 0$ implies
$\lambda:\Delta^n\rightarrow\reals$ is monotone decreasing. This
proves $\phi$ is monotone. It is Scott continuous because
it is Euclidean continuous and suprema in the Bayesian order
are pointwise Euclidean limits. Finally,
\[\lambda(\phi x)=I_f(\lambda x)\]                        
and so by induction
\[\lambda(\phi^n x)=I_f^n(\lambda x).\]
Then its least fixed point $\fix(\phi)$ must satisfy
\[\lambda(\fix(\phi))=\lambda(\bigsqcup\phi^n(\bot))=\lambda(\lim_{n\rightarrow\infty}\phi^n(\bot))=
\lim_{n\rightarrow\infty}\lambda(\phi^n(\bot))=\lim_{n\rightarrow\infty}I_f^n(0)=\lambda\]
which gives
\[\fix(\phi)=\phi(\fix(\phi))=(e^{\lambda a_1},\ldots,e^{\lambda a_n})\cdot\frac{1}{\sum_{i=1}^n e^{\lambda a_i}}\]
the maximum entropy state.
\endproof\newline

The map $\phi$ is not monotone
with respect to majorization $(\Lambda^n,\sqsubseteq)$. To see
why, take a problem where $\lambda = 0$,
then $I_f(0)=0$, which means $\lambda x\leq 0$. 
Then $\phi:\Lambda^n\rightarrow\Lambda^n$. Let $x=(1/2,2/5,1/10)$,
$y=(1/2,1/2,0)$. Then $x\sqsubseteq y$ in majorization.
Because $\lambda x<0$, $I_f(\lambda x)<0$, so 
$\phi(x)\neq\bot$. However, $\phi(y)=\bot$,
which means $\phi(x)\not\sqsubseteq\phi(y)$ in majorization.
It is not immediately clear
whether $\phi$ is monotone in the implicative
order~\cite{cont}. Notice though that Lemma~\ref{lambdalemma} is
also valid for the implicative order.

\section{A conspiracy theory}

The Bayesian projections $(p_i)$ used to define
the Bayesian order relate to entropy in a special way:
\[\sigma x = (1-x_k)\sigma p_k(x) + \sigma(x_k,1-x_k)\]
for any $k$ with $x_k\neq 1$. This property 
might imply Shannon's additivity
property or `the recursion axiom' so that any function 
satisfying the equation above
and the two other usual axioms has to be entropy to within a constant.
This equation looks like it almost means something.
\em If we knew what it meant\em, and we could also use it to
establish a plausible link to the Bayesian order, 
then we might try to prove this uniqueness.

\section{Etc.}

The things in
this paper which are original (to the best
of our knowledge) are that $I_f$ always
iterates to $\lambda$ and that it 
can be used to define the Scott continuous
$\phi$ whose least fixed point is the
maximum entropy state.

The
operator $\phi$ might have a meaningful 
interpretation: we begin with $\bot$, and then with
each iteration probabilities are adjusted
based on the information $(\langle a|x\rangle=E)$ 
we have, until the limit gives us
just the right state. There should be a logic that captures
the type of inference provided by
the maximum entropy principle: states are propositions.
Perhaps the logic we are looking
for treats observables 
as \em incomplete descriptions \em of propositions:
\begin{enumerate}

\item[(i)] Maybe 
the logic
has sequents of the form $I\rightarrow q$
where $q$ is a proposition and $I$ is
information which partially describes
a proposition. In this logic, it should
be a theorem that $a,E\rightarrow\fix(\phi)$.
Is it possible to extract a `proof' 
of this theorem from $(\phi^n\bot)$? Or
is the theorem $\fix(\phi)\rightarrow\langle a\rangle=E$?
What is the logic of statistical mechanics?

\item[(ii)] The maximum entropy principle
and some of its variants
might all have some underlying qualitative component. Maybe
it is possible to explain how from a certain kind of 
logic one can extract a statistical inference method.
Maybe one has a choice about how to write 
the maximum entropy principle: either in the language of
expectations, entropy and optimization techniques,
or as a logic that will probably annoy a lot of people.

\item[(iii)] Maybe $\phi$ has an informatic derivative at $\fix(\phi)$. 
\end{enumerate}

\end{document}